\documentclass[12pt,a4paper,reqno]{amsart}

\usepackage{geometry}
\title[Decomposition of elliptic multiple zeta values and iterated Eisenstein integrals]{Decomposition of elliptic multiple zeta values}
\author{Nils Matthes}

\usepackage{amssymb,amsmath,amscd,amsthm,mathrsfs}

\usepackage[bookmarksopen]{hyperref}

\usepackage{color}
\usepackage{shuffle}

\usepackage{mathtools}
\mathtoolsset{showonlyrefs}

\date{}
\address{Fachbereich Mathematik (AZ)\\ Universit\"at Hamburg\\ Bundesstrasse 55\\ D-20146 Hamburg}
\email{nils.matthes@uni-hamburg.de}
\subjclass[2010]{11M32, (11F67)}
\keywords{Modular symbols, elliptic associators, elliptic polylogarithms}

\makeindex

\newif\ifnote 
\notetrue

\allowdisplaybreaks

\theoremstyle{definition}
\newtheorem{dfn}{Definition}[section]
\newtheorem{rmk}[dfn]{Remark}
\newtheorem{exmp}[dfn]{Example}

\theoremstyle{plain}
\newtheorem{prop}[dfn]{Proposition}

\newtheorem{thm}[dfn]{Theorem}

\newenvironment{prf}{\begin{proof}[{\it Proof: \nopunct}]}{\end{proof}}

\numberwithin{equation}{section}

\def\dd{\mathrm{d}}

\def\bC{\mathbb C}
\def\bZ{\mathbb Z}
\def\bQ{\mathbb Q}

\def\bH{\mathbb H}

\def\bR{\mathbb R}

\def\cE{\mathcal E}

\def\cL{\mathcal L}

\def\cO{\mathcal O}

\def\cZ{\mathcal Z}

\def\fm{\mathfrak{m}}

\def\fu{\mathfrak{u}}

\def\eZ{\mathcal{EZ}}

\DeclareMathOperator{\ad}{ad}

\DeclareMathOperator{\Der}{Der}

\DeclareMathOperator{\im}{Im}
\DeclareMathOperator{\Lie}{Lie}
\DeclareMathOperator{\SL}{SL}

\DeclareMathOperator{\spn}{Span}

\def\be{\mathbf{e}}

\def\ul{\underline}

\def\bfE{\mathbf{E}}
\def\bfF{\mathbf{F}}
\begin{document}

\begin{abstract}
	\noindent
	We describe a decomposition algorithm for elliptic multiple zeta values, which amounts to the construction of an injective map $\psi$ from the algebra of elliptic multiple zeta values to a space of iterated Eisenstein integrals. We give many examples of this decomposition, and conclude with a short discussion about the image of $\psi$. It turns out that the failure of surjectivity of $\psi$ is in some sense governed by period polynomials of modular forms.
\end{abstract}

\title[Decomposition of elliptic multiple zeta values]{Decomposition of elliptic multiple zeta values and iterated Eisenstein integrals}

\maketitle

\section{Introduction} \label{sec:1}

The purpose of this paper is to describe a decomposition of elliptic multiple zeta values into linear combinations of iterated Eisenstein integrals, which clarifies the algebraic structure of elliptic multiple zeta values. This point of view was first taken up in \cite{BMS}, and the present paper summarizes and extends some of its main results.
\subsection{Elliptic multiple zeta values} \label{ssec:1.1}
The notion of elliptic multiple zeta value first appeared explicitly in \cite{Eemzv} under the name ``analogues elliptiques des nombres multiz\'etas''.
Elliptic multiple zeta values come in two closely related versions, namely A-elliptic and B-elliptic multiple zeta values, each corresponding to one of the two canonical homology cycles $\alpha$ and $\beta$ on a once-punctured elliptic curve.
Elliptic multiple zeta values are linked to a variety of other subjects, such as multiple elliptic polylogarithms \cite{BL,LR}, elliptic braid groups and elliptic associators \cite{CEE,Eass}, as well as mixed elliptic motives \cite{HM}. They also occur in amplitude computations in string theory \cite{BMMS}.
%A recurring theme in their study is the analogy with the classical multiple zeta values, a first glimpse of which can already be seen from their definition: just as multiple zeta values can be defined as iterated integrals on the punctured projective line $\bP_{\bC}^1 \setminus \{0,1,\infty\}$, elliptic multiple zeta values are defined by iterated integrals on a once-punctured elliptic curve $\bC/(\bZ+\bZ\tau) \setminus \{0\}$. Also, the notions of weight and depth for multiple zeta values have analogue for elliptic multiple zeta values, namely the weight and the length.

Elliptic multiple zeta values are known to satisfy very many $\bQ$-linear relations, which are studied in \cite{BMMS,BMS,edzv}. Understanding the entirety of all such relations is a delicate problem, which, despite some advances \cite{edzv}, is not yet fully understood.

The main theme of this paper is that the study of relations becomes somewhat simpler when one rewrites elliptic multiple zeta values as iterated integrals on the upper half-plane. More precisely, every elliptic multiple zeta value can be decomposed uniquely as a linear combination of \textit{iterated Eisenstein integrals} \cite{MMV,Man}. The gain of this representation is that the set of all iterated Eisenstein integrals is linearly independent over $\bC$ \cite{LMS}, thus finding relations between elliptic multiple zeta values reduces to solving linear systems of equations. In fact, this procedure was used in \cite{edzv} to prove optimal lower bounds for spaces of elliptic double zeta values.
\subsection{Analogy with decomposition of motivic multiple zeta values} \label{ssec:1.2}
The decomposition of elliptic multiple zeta values into iterated Eisenstein integrals is in many ways reminiscent of the decomposition algorithm for motivic multiple zeta values into polynomials in non-commutative variables $f_3,f_5,f_7,\ldots$, the so-called ``$f$-alphabet'', which amounts to an isomorphism  \cite{mtm,Brown:decomp}
\begin{equation}
\phi: \cZ^{\mathfrak{m}} \stackrel{\cong}\longrightarrow T(\bfF)^{\vee} \otimes_{\bQ} \bQ[f_2].
\end{equation}
Here $\cZ^{\mathfrak m}$ is the $\bQ$-algebra of motivic multiple zeta values $\zeta^{\fm}(k_1,\ldots,k_n)$, $T(\bfF)^{\vee}$ denotes the graded dual of the tensor algebra on the set $\bfF=\{f_{2k+1} \, \vert \, k\geq 1\}$, and $f_2$ is an additional variable, which commutes with all $f_{2k+1}$ and corresponds to $\zeta^{\fm}(2)$. The main ingredient for the construction of $\phi$ is the motivic coaction for motivic multiple zeta values \cite{mtm,GonGalSym}. In addition, one needs to choose a set of free algebra generators for $\cZ^{\mathfrak m}$, and therefore the construction of $\phi$ is not canonical.

The analogous situation for elliptic multiple zeta values is similar, but technically simpler. Denoting by $\eZ^{\rm A}$ the $\bQ$-algebra of A-elliptic multiple zeta values, there is an injection \cite{BMS,thesis}
\begin{equation}
\psi^{\rm A}: \eZ^{\rm A} \hookrightarrow T(\bfE)^{\vee} \otimes_{\bQ} \cZ[2\pi i],
\end{equation}
where $\bfE=\{\be_0,\be_2,\be_4\ldots\}$ and $\cZ$ is the $\bQ$-algebra of multiple zeta values.\footnote{One also has an injection $\psi^{\rm B}: \eZ^{\rm B} \hookrightarrow T(\bfE)^{\vee} \otimes_{\bQ} \cZ[2\pi i]$, where $\eZ^{\rm B}$ is the algebra of B-elliptic multiple zeta values \cite{thesis}.} Here, the variable $\be_{2k}$ should be thought of as corresponding to the Eisenstein series $E_{2k}(\tau)$ for $\SL_2(\bZ)$ (where $E_0(\tau):=-1$). The key to the construction of $\psi^{\rm A}$ is the differential equation for elliptic multiple zeta values, found by Enriquez \cite{Eass,Eemzv}. In fact, we argue that this differential equation can be seen as an elliptic analogue of the motivic coaction. In contrast to the $\phi$-map for motivic multiple zeta values, the construction of $\psi^{\rm A}$ is completely canonical and does not depend on any initial choices. However, unlike $\phi$, the morphism $\psi^{\rm A}$ is not an isomorphism: The failure of surjectivity is related to a certain Lie algebra of derivations, and ultimately to the existence of modular forms for $\SL_2(\bZ)$ \cite{BMMS,Pol}. A precise description of the image of $\psi^{\rm A}$ will be given in a joint work with Lochak and Schneps \cite{LMS}.
\subsection{Overview of the article} \label{ssec:1.3}
In Section \ref{sec:2}, we give a brief introduction to iterated Eisenstein integrals, focusing on their algebraic structure. Section \ref{sec:3} contains the definition of elliptic multiple zeta values and also a short discussion of their differential equation. The consequences of this differential equation are then studied in the remaining sections: In Section \ref{sec:4}, which essentially follows \cite{BMS}, we construct the map $\psi^{\rm A}$ and give many concrete examples. Then, in Section \ref{sec:5}, we turn our attention towards describing the image of $\psi^{\rm A}$ by relating it to the aforementioned Lie algebra of derivations.
\subsection{Acknowledgments}
This article was written on the occasion of the conference \textit{Various aspects of multiple zeta values}, held at the \textit{Research Institute for Mathematical Sciences} (RIMS) in Kyoto during July 2016. It is my pleasure to thank the organizer Hidekazu Furusho for the invitation to talk there. Also, many thanks to the mathematical department of Nagoya University and RIMS for hospitality.

\section{Iterated Eisenstein integrals} \label{sec:2}

Iterated Eisenstein integrals are a special case of iterated Shimura integrals \cite{Man}, which were introduced by Manin to study the rational homotopy theory of modular curves. More recently, the theory of iterated Shimura integrals has been thoroughly revisited and extended by Brown \cite{MMV}.
In the context of this paper, iterated Eisenstein integrals will be the basic building blocks of elliptic multiple zeta values.

\subsection{Generalities} \label{ssec:2.1}
We begin by recalling the general notion of an iterated integral, due to Chen \cite{Ch}. 
Let $M$ be a complex manifold. Given a collection of smooth differential one-forms $\omega_1,\ldots,\omega_n \in \Omega^1(M)$ and a piecewise smooth path $\gamma: [0,1] \rightarrow M$, one defines the iterated integral
\begin{equation}
\int_{\gamma}\omega_1\ldots\omega_n:=\int_{0\leq t_1\leq \ldots \leq t_n\leq 1}\gamma^*(\omega_1)(t_1)\ldots\gamma^*(\omega_n)(t_n) \in \bC,
\end{equation}
where $\gamma^*(\omega)(t) \in \Omega^1([0,1])$ denotes the pull-back of $\omega$ along $\gamma$, and $t$ is the natural coordinate on $[0,1]$. If $n=0$, then $\int_{\gamma}:=1$ (the empty iterated integral).

General properties of iterated integrals include the shuffle product formula
\begin{equation} \label{eqn:shuffle}
\int_{\gamma}\omega_1\ldots\omega_r\int_{\gamma}\omega_{r+1}\ldots\omega_{r+s}=\sum_{\sigma \in \Sigma_{r,s}}\int_{\gamma}\omega_{\sigma(1)}\ldots\omega_{\sigma(r+s)},
\end{equation}
where $\Sigma_{r,s} \subset \Sigma_{r+s}$ denotes the set of $(r,s)$-shuffle, i.e. the set of all permutations $\sigma$ of the set $\{1,\ldots,r+s\}$ such that $\sigma^{-1}$ is strictly increasing on both $\{1,\ldots,r\}$ and on $\{r+1,\ldots,r+s\}$.
Also, iterated integrals satisfy the differential equation
\begin{equation} \label{eqn:diffeq}
\frac{\dd}{\dd t}\Bigr|_{t=a}\int_{\gamma_t}\omega_1\ldots\omega_n=-\langle\omega_1,\gamma'(a)\rangle \int_{\gamma_a}\omega_2\ldots\omega_n,
\end{equation}
where $\langle \cdot,\cdot \rangle$ is the natural pairing and for $a \in [0,1]$, we denote by $\gamma_a: [0,1] \rightarrow M$ the path $t \mapsto \gamma(t+(1-t)a)$. For more properties of iterated integrals, we refer to \cite{Colombia,Hai}.
\subsection{Iterated integrals on the upper half-plane}
In the definition of iterated integrals, we will be mainly interested in the case where $M$ is the upper half-plane $\bH=\{z \in \bC \, \vert \, \im(z)>0\}$. In this case, if the differential one-forms $\omega_1,\ldots,\omega_n$ are holomorphic, the value of the iterated integral $\int_{\gamma}\omega_1\ldots\omega_n$ depends only on the start and end point of $\gamma$ (this holds more generally on every one-dimensional and simply connected complex manifold). Hence, given two points $a,b \in \bH$ and $\omega_1,\ldots,\omega_n$ as above, we may write $\int_a^b\omega_1\ldots\omega_n$ without ambiguity.

One can also define iterated integrals along a path between a point $\tau \in \bH$ and the cusp $i\infty$, provided the differential forms $\omega_i$ have at most simple poles at $i\infty$. This uses Deligne's tangential base points (cf. \cite{Del}, \S 15), and is worked out in detail in the case of iterated integrals on $\bH$ in \cite{MMV}, Section 4. In the sequel, we use the conventions and notation from \cite{MMV}, in particular, all our integrals are regularized with respect to the tangent vector $\overrightarrow{1}_{\infty}$ at $i\infty$.

The iterated integrals on $\bH$ we are interested in are the \textit{iterated Eisenstein integrals}
\begin{align} \label{eqn:iei}
\cE(2k_1,\ldots,2k_n;\tau):=(2\pi i)^n\int_{\tau}^{i\infty}E_{2k_1}(\tau_1)\dd\tau_1\ldots E_{2k_n}(\tau_n)\dd\tau_n,
\end{align}
where $k_1,\ldots,k_n \geq 0$ and $\tau \in \bH$. Here, $E_0(\tau):=-1$ and for $k \geq 1$, $E_{2k}(\tau)$ denotes the Hecke-normalized Eisenstein series\footnote{The case $k=1$ requires Eisenstein summation
\[
\sum_{(m,n) \in \bZ^2}a_{m,n}:=\lim_{N \to \infty}\lim_{M \to \infty}\sum_{n=-N}^N\sum_{m=-M}^M a_{m,n}.
\]
}
\begin{align}\label{eqn:Eisen}
E_{2k}(\tau)=\frac{(2k-1)!}{2(2\pi i)^{2k}}\sum_{(m,n) \in \bZ^2 \setminus \{(0,0)\}}\frac{1}{(m+n\tau)^{2k}}=-\frac{B_{2k}}{4k}+\sum_{n \geq 1}\sigma_{2k-1}(n)q^n,
\end{align}
%\begin{equation} \label{eqn:Eisen}
%E_{2k}(\tau)=-\frac{B_{2k}}{4k}+\sum_{n \geq 1}\sigma_{2k-1}(n)q^n, \quad q=e^{2\pi i\tau},
%\end{equation}
where $q=e^{2\pi i\tau}$, $B_m$ denotes the $m$-th Bernoulli number, defined by $\frac{t}{e^t-1}=\sum_{m\geq 0}B_m\frac{t^m}{m!}$, and $\sigma_m(n):=\sum_{d \vert n}d^m$ is the $m$-th divisor function.
For notational convenience, we extend the definition of the Eisenstein series to all non-negative integers by setting $E_k(\tau)=0$, if $k \geq 1$ is odd. In particular, $\cE(k_1,\ldots,k_n;\tau)=0$, if one of the $k_i$ is odd.

The iterated Eisenstein integrals $\cE(k_1,\ldots,k_n;\tau)$ are holomorphic functions of $\tau$ and the analogue of the differential equation \eqref{eqn:diffeq} is (cf. \cite{MMV}, Proposition 4.7)
\begin{equation} \label{eqn:ieidiffeq}
\frac{1}{2\pi i}\frac{\dd}{\dd \tau}\Bigr|_{\tau=\rho}\cE(k_1,\ldots,k_n;\tau)=-E_{k_1}(\rho)\cE(k_2,\ldots,k_n;\rho).
\end{equation}

\subsection{The algebra of iterated Eisenstein integrals} \label{ssec:2.2}
Let $\bQ\langle\cE\rangle \subset \cO(\bH)$ be the $\bQ$-vector subspace spanned by the iterated Eisenstein integrals (where $\cO(\bH)$ is the $\bC$-algebra of holomorphic functions on $\bH$). By \eqref{eqn:shuffle}, we have the shuffle product formula (cf. e.g. \cite{MMV}, Proposition 4.7)
\begin{equation} \label{eqn:shuffleeis}
\cE(k_1,\ldots,k_r;\tau)\cE(k_{r+1},\ldots,k_{r+s};\tau)=\sum_{\sigma \in \Sigma_{r,s}}\cE(k_{\sigma(1)},\ldots,k_{\sigma(r+s)};\tau).
\end{equation}
In particular, $\bQ\langle\cE\rangle$ is a $\bQ$-algebra. In order to describe $\bQ\langle\cE\rangle$ in more detail, let $\mathbf{E}:=\{\be_{0},\be_2,\be_4,\ldots\}$ be a set of variables indexed by the non-negative even integers, and let $T(\mathbf{E})$ be the tensor $\bQ$-algebra, which is graded by giving the variables $\be_{2k}$ degree one. In fact, $T(\mathbf{E})$ has the natural structure of a graded Hopf algebra: its coproduct $\Delta: T(\mathbf{E}) \rightarrow T(\mathbf{E}) \otimes_{\bQ} T(\mathbf{E})$ is the unique coproduct such that all the $\be_{2k}$ are primitive, i.e. $\Delta(\be_{2k})=\be_{2k} \otimes_{\bQ}1+1\otimes_{\bQ}\be_{2k}$ for all $k \geq 0$, and its antipode is the unique anti-homomorphism sending $\be_{2k} \mapsto -\be_{2k}$. We denote by $T(\mathbf{E})^{\vee}$ the graded dual of $T(\mathbf{E})$, which is the Hopf algebra dual of $T(\mathbf{E})$. Its product is the shuffle product
\begin{align}
\shuffle: T(\mathbf{E})^{\vee} \otimes_{\bQ} T(\mathbf{E})^{\vee} &\rightarrow T(\mathbf{E})^{\vee}\notag\\
\be^{\vee}_{2k_1}\ldots\be^{\vee}_{2k_r} \otimes_{\bQ} \be^{\vee}_{2k_{r+1}}\ldots\be^{\vee}_{2k_{r+s}} &\mapsto \sum_{\sigma \in \Sigma_{r,s}}\be^{\vee}_{2k_{\sigma(1)}}\ldots \be^{\vee}_{2k_{\sigma(r+s)}},
\end{align}
and its coproduct is given by deconcatenation
\begin{align}
\Delta^{\vee}: T(\mathbf{E})^{\vee} &\rightarrow T(\mathbf{E})^{\vee} \otimes_{\bQ} T(\mathbf{E})^{\vee} \notag\\
\be_{2k_1}^{\vee}\ldots\be_{2k_n}^{\vee} &\mapsto \sum_{i=0}^n\be_{2k_1}^{\vee}\ldots\be_{2k_i}^{\vee} \otimes_{\bQ}\be_{2k_{i+1}}^{\vee}\ldots\be_{2k_n}^{\vee}.
\end{align}
Given a multi-index $\underline{2k}=(2k_1,\ldots,2k_n) \in (2\bZ_{\geq 0})^n$, we will frequently write $\be_{\underline{2k}}$ instead of $\be_{2k_1}^{\vee}\ldots\be_{2k_n}^{\vee}$.

The following theorem is a consequence of $\bC$-linear independence of iterated Eisenstein integrals \cite{LMS}. It shows in particular that $\bQ\langle\cE\rangle$ is a graded Hopf algebra in a natural way.
\begin{thm} \label{thm:ieiisom}
For any $\bQ$-subalgebra $K \subset \bC$, the $K$-linear morphism
\begin{align} \label{eqn:ieiisom}
\psi: \bQ\langle\cE\rangle \otimes_{\bQ}K &\rightarrow T(\mathbf{E})^{\vee} \otimes_{\bQ}K,\notag\\
\cE(2k_1,\ldots,2k_n;\tau) &\mapsto \be_{2k_1}^{\vee}\ldots\be_{2k_n}^{\vee}
\end{align}
is a well-defined isomorphism of $K$-algebras.
\end{thm}
In particular, the only algebraic relations between iterated Eisenstein integrals are given by \eqref{eqn:shuffleeis}.
\section{Elliptic multiple zeta values} \label{sec:3}
In this section, we recall the definition of elliptic multiple zeta values \cite{Eemzv} (see also \cite{BMMS,BMS,edzv}). An important role in the definition is played by a certain Jacobi form in two variables, whose study dates back to Eisenstein and Kronecker \cite{W}.
\subsection{Differential forms on a once-punctured elliptic curve} \label{ssec:3.1}
For a point $\tau \in \bH$, we will denote by $E_{\tau}^{\times}:=\bC/(\bZ+\bZ\tau) \setminus \{0\}$ the associated once-punctured complex elliptic curve, with its canonical coordinate $\xi=s+r\tau$, where $r,s \in \bR$. In \cite{BL}, Brown and Levin have introduced the following differential one-form
\begin{equation} \label{eqn:BL}
\Omega_{\tau}(\xi,\alpha)=e^{2\pi ir\alpha}\frac{\theta_{\tau}'(0)\theta_{\tau}(\xi+\alpha)}{\theta_{\tau}(\xi)\theta_{\tau}(\alpha)}\dd\xi,
\end{equation}
which is a variant of the Kronecker-Eisenstein series $F_{\tau}(\xi,\alpha)=\frac{\theta_{\tau}'(0)\theta_{\tau}(\xi+\alpha)}{\theta_{\tau}(\xi)\theta_{\tau}(\alpha)}$ \cite{W,Zag}. Here,
\begin{equation}
\theta_{\tau}(\xi)=\sum_{n\in \bZ}(-1)^ne^{2\pi i\xi\left(n+\frac 12\right)}e^{\pi i\tau\left(n+\frac 12\right)^2}
\end{equation}
is the odd Jacobi theta function.
%Note that the function $F_{\tau}(\xi,\alpha)=\frac{\theta_{\tau}'(0)\theta_{\tau}(\xi+\alpha)}{\theta_{\tau}(\xi)\theta_{\tau}(\alpha)}$ is meromorphic with simple poles along the hyperplanes $\xi=0$ and $\alpha=0$ (cf. \cite{Zag}, Theorem in Section 3). In particular, it has a formal expansion at $\alpha=0$, and the same is true for $\Omega_{\tau}(\xi,\alpha)$. This expansion gives rise to a family $\{\omega^{(k)}\}$ of smooth differential one-forms on $E_{\tau}^{\times}$
As explained in \cite{BL}, Section 3.5, $\Omega_{\tau}(\xi,\alpha)$ is invariant under lattice translations $\xi \mapsto \xi+m+n\tau$ for $m,n \in \bZ$, and has a formal expansion in $\alpha$
\begin{equation}
\Omega_{\tau}(\xi,\alpha)=\sum_{k\geq 0}\omega^{(k)}\alpha^{k-1},
\end{equation}
where every $\omega^{(k)}$ is a smooth differential one-form on $E_{\tau}^{\times}$.
\subsection{Definition and first examples of elliptic multiple zeta values} \label{ssec:3.2}
Elliptic multiple zeta values will be defined as regularized iterated integrals of the forms $\omega^{(k)}$ along paths on $E_{\tau}^{\times}$. There are two natural such choices, namely the images $\alpha$ and $\beta$ of the (open) straight line paths from $0$ to $1$ resp. from $0$ to $\tau$ under the projection $\bC \setminus (\bZ+\bZ\tau) \rightarrow E_{\tau}^{\times}$. Corresponding to the two natural paths $\alpha$ and $\beta$ on the once-punctured elliptic curve $E_{\tau}^{\times}$, there are two types of elliptic multiple zeta values, namely A-elliptic and B-elliptic multiple zeta values, which are related to one another by a certain modular transformation formula (cf. \cite{Eemzv}, Section 2.5). For simplicity, we will consider in this paper only the A-elliptic multiple zeta values.
\begin{dfn} \label{dfn:eMZV}
For integers $k_1,\ldots,k_n \geq 0$, define the \textit{A-elliptic multiple zeta value} $I^{\rm A}(k_1,\ldots,k_r;\tau)$ to be the regularized\footnote{See \cite{edzv}, Definition 2.1 for the details, which employs Deligne's regularization prescription using tangential base points (\cite{Del}, \S 15).} iterated integral
\begin{equation}
I^{\rm A}(k_1,\ldots,k_n;\tau)=(2\pi i)^{-(k_1+\ldots+k_n-n)}\int_{\alpha}\omega^{(k_1)}\ldots\omega^{(k_n)}.
\end{equation}
The \textit{length} of $I^{\rm A}(k_1,\ldots,k_n;\tau)$ is defined to be $n$.
\end{dfn}
\begin{rmk}
The original reference for elliptic multiple zeta values is \cite{Eemzv}, with additional references being \cite{BMMS,BMS,edzv}. Note that the pre-factor $(2\pi i)^{-(k_1+\ldots+k_n-n)}$ is not included in the original definition of A-elliptic multiple zeta values. In the context of this paper, introducing this factor has the benefit of removing many cumbersome powers of $2\pi i$ from the formulas, which will make the algebraic structure of A-elliptic multiple zeta values more transparent.
\end{rmk}
As functions of $\tau$, A-elliptic multiple zeta values are holomorphic on the upper half-plane. In fact, more is true.
\begin{prop}[\cite{Eemzv}, Proposition 5.3] \label{prop:Fourier}
Every A-elliptic multiple zeta value has a convergent Fourier expansion
\begin{equation}
\sum_{m \geq 0}a_mq^m, \quad q=e^{2\pi i\tau},
\end{equation}
such that $a_m \in \cZ[2\pi i]$, where $\cZ$ denotes the $\bQ$-algebra of multiple zeta values.
\end{prop}
\begin{exmp} \label{exmp:lengthone}
In length one, we have
\begin{equation} \label{eqn:lengthone}
I^{\rm A}(k;\tau)=\begin{cases}\frac{2\pi iB_k}{k!} & \mbox{if $k$ is even,}\\ 0 &\mbox{if $k$ is odd,}\end{cases}
\end{equation}
where $B_k$ denotes the $k$-th Bernoulli number. This is straightforward to verify using the definition $I^{\rm A}(k;\tau)=\int_{\alpha}\omega^{(k)}$ and the Fourier expansion of the Kronecker-Eisenstein series $F_{\tau}(\xi,\alpha)$ (cf. \cite{Zag}, Theorem 3).
\end{exmp}

\subsection{Differential equation and constant term} \label{ssec:3.3}
In \cite{Eemzv}, Th\'eor\`eme 3.3, Enriquez has found a differential equation for elliptic multiple zeta values, viewed as holomorphic functions in the coordinate $\tau$ on $\bH$. This differential equation is recursive for the length of the elliptic multiple zeta values, and can be expressed using Eisenstein series.\footnote{To be precise, the result in \cite{Eemzv} is expressed in terms of the (not normalized) Eisenstein series $G_{2k}(\tau)=\frac{2(2\pi i)^{2k}}{(2k-1)!}E_{2k}(\tau)$.}
\begin{thm}[Enriquez] \label{thm:diffeq}
We have
\small{
\begin{align}
&\frac{1}{2\pi i}\frac{\dd}{\dd \tau} I^{\rm A}(k_1,\ldots,k_n;\tau) =  \alpha_{k_1+1} E_{k_1+1}(\tau) I^{\rm A}(k_2,\ldots,k_n;\tau) - \alpha_{k_n+1} E_{k_n+1}(\tau) I^{\rm A}(k_1,\ldots,k_{n-1};\tau)\notag \\
&+ \sum_{i=2}^n \Bigg\{ (-1)^{k_i} \alpha_{k_{i-1}+k_i+1} E_{k_{i-1}+k_i+1}(\tau) I^{\rm A}(k_1,\ldots,k_{i-2},0,k_{i+1},\ldots,k_n;\tau) \label{eqn:tauder} \\
&-  \sum_{m=0}^{k_{i-1}+1} \binom{k_i+m-1}{m}  \alpha_{k_{i-1}-m+1}E_{k_{i-1}-m+1}(\tau) I^{\rm A}(k_1,\ldots,k_{i-2},m+k_i,k_{i+1},\ldots,k_n;\tau) \notag \\
&+ \sum_{m=0}^{k_{i}+1} \binom{k_{i-1}+m-1}{m}  \alpha_{k_i-m+1}E_{k_i-m+1}(\tau) I^{\rm A}(k_1,\ldots,k_{i-2},m+k_{i-1},k_{i+1},\ldots,k_n;\tau) \Bigg\},
\notag 
\end{align}
}
where $\alpha_n:=\begin{cases}-1 & \mbox{if $n=0$,} \\ 0 & \mbox{if $n=1$,} \\ \frac{2}{(n-2)!} & \mbox{if $n\geq 2$.}\end{cases}$
\end{thm}
Given $I^{\rm A}(k_1,\ldots,k_n;\tau)$ with Fourier expansion $\sum_{m \geq 0}a_mq^m$, we have 
\begin{equation}
\frac{1}{2\pi i}\frac{\dd}{\dd\tau}I^{\rm A}(k_1,\ldots,k_n;\tau)=\sum_{m\geq 1}ma_mq^m.
\end{equation}
Thus, \eqref{eqn:tauder} gives a recursive formula for the Fourier coefficients $a_m$ for $m \geq 1$. On the other hand, the constant term $a_0$ in the Fourier expansion is given by $\lim_{\tau \to i\infty}I^{\rm A}(k_1,\ldots,k_n;\tau)$.

In order to retrieve the constant terms of A-elliptic multiple zeta values in a systematic way, we consider the generating series of A-elliptic multiple zeta values
\begin{equation} \label{eqn:Aassoc}
\underline{A}(\tau):=\sum_{n\geq 0} (-1)^n\sum_{k_1,\ldots,k_n \geq 0}I^{\rm A}(k_1,\ldots,k_n;\tau)\ad^{k_n}(a)(b)\ldots \ad^{k_1}(a)(b) \in \bC\langle\!\langle a,b\rangle\!\rangle.
\end{equation}
Here, $\bC\langle\!\langle a,b\rangle\!\rangle$ is the $\bC$-algebra of formal power series in the non-commuting variables $a$ and $b$, and $\ad^k(a)$ denotes the $k$-fold iterate of the adjoint action $\ad(a)(p)=ap-pa$ on $\bC\langle\!\langle a,b\rangle\!\rangle$. The series $\underline{A}(\tau)$ is related to the series $A(\tau) \in \bC\langle\!\langle a,b\rangle\!\rangle$ occurring as a datum of Enriquez's elliptic KZB associator \cite{Eass} by $\underline{A}(\tau)=e^{-\pi i[a,b]}A(\tau)$.\footnote{To be precise, Enriquez writes the elliptic KZB associator in variables $x$, $y$, which are related to the variables $a$, $b$ introduced here by $a=2\pi ix$, $b=(2\pi i)^{-1}y$. This also slightly changes the appearance, though not the essential content, of several results concerning $\underline{A}(\tau)$ such as Theorems \ref{thm:constterm} and \ref{thm:diffeqass}.}
\begin{thm}[\cite{Eass}, Proposition 6.3] \label{thm:constterm}
The limit $\lim_{\tau \to i\infty}\underline{A}(\tau)$ exists and we have
\begin{equation} \label{eqn:constterm}
\lim_{\tau \to i\infty}\underline{A}(\tau)=e^{\pi it}\Phi(\tilde{y},t)e^{2\pi i\tilde{y}}\Phi(\tilde{y},t)^{-1},
\end{equation}
where $t=-[a,b]$ and $\tilde{y}=-\frac{\ad(a)}{e^{\ad(a)}-1}(b)$ and $\Phi$ is the Drinfeld associator \cite{Dr,FurStab}.
\end{thm}
%Extracting a constant term for a specific A-elliptic multiple zeta value from this theorem is cumbersome to do by hand, but it can be done using an implementation in a suitable CAS (in our case, Mathematica\textsuperscript{\textregistered}).
The coefficients of the Drinfeld associator $\Phi$ are given by $\bQ$-linear combinations of multiple zeta values; see \cite{FurStab}, Proposition 3.2.3 for an explicit formula. Comparing coefficients on both sides of \eqref{eqn:constterm}, one therefore obtains a formula for the constant term $\lim_{\tau \to i\infty}I^{\rm A}(k_1,\ldots,k_n;\tau)$ of $I^{\rm A}(k_1,\ldots,k_n;\tau)$ in terms of $\bQ[2\pi i]$-linear combinations of multiple zeta values, which is however rather cumbersome to write down in practice (cf. \cite{BMS}, Section 2.3.1 for some examples).

It turns out that the differential equation for A-elliptic multiple zeta values (i.e. Theorem \ref{thm:diffeq}) can also be expressed using the generating series $\underline{A}(\tau)$. We will come back to this in Section \ref{sec:5}.
\section{Decomposition of elliptic multiple zeta values} \label{sec:4}

We will show how the results of the last section can be used to rewrite A-elliptic multiple zeta values as iterated Eisenstein integrals. This has the crucial advantage that, by Theorem \ref{thm:ieiisom}, the algebraic relations satisfied by iterated Eisenstein integrals are completely under control.
\subsection{The decomposition map} \label{ssec:4.2}
The starting point is the interpretation of the differential equation \eqref{eqn:tauder} as a statement about the algebraic structure of A-elliptic multiple zeta values. Let $\eZ^{\rm A}$ be the $\bQ$-vector space spanned by the A-elliptic multiple zeta values
\begin{equation}
\eZ^{\rm A}:=\spn_{\bQ}\{I^{\rm A}(k_1,\ldots,k_n;\tau) \, \vert \, n\geq 0, \, k_i \geq 0 \}.
\end{equation}
By the shuffle product formula for iterated integrals \eqref{eqn:shuffle}, $\eZ^{\rm A}$ is a $\bQ$-algebra.
\begin{prop} \label{prop:psiA}
There is a natural embedding of $\bQ$-algebras
\begin{equation}
\psi^{\rm A}: \eZ^{\rm A} \hookrightarrow T(\mathbf{E})^{\vee} \otimes_{\bQ} \cZ[2\pi i].
\end{equation}
%which is the restriction of the isomorphism $\psi: \bQ\langle\cE\rangle \otimes_{\bQ} \cZ[2\pi i] \rightarrow T(\mathbf{E})^{\vee} \otimes_{\bQ} \cZ[2\pi i]$ of Theorem \ref{thm:ieiisom}.
\end{prop} 
\begin{prf}
We first claim that every A-elliptic multiple zeta value can be written as a $\cZ[2\pi i]$-linear combination of iterated Eisenstein integrals \eqref{eqn:iei}. By Example \ref{exmp:lengthone}, we have $I^{\rm A}(2k;\tau)=\frac{2\pi iB_{2k}}{(2k)!}$ and $I^{\rm A}(2k+1;\tau)=0$, which are $\cZ[2\pi i]$-linear combinations of the empty iterated Eisenstein integral $\cE(\emptyset;\tau)=1$, hence the claim is true for A-elliptic multiple zeta values of length one. Now assume the claim for all A-elliptic multiple zeta values up to and including length $n-1$. By the differential equation \eqref{eqn:tauder}, we know that $\frac{1}{2\pi i}\frac{\dd}{\dd \tau}I^{\rm A}(k_1,\ldots,k_n;\tau)$ is a $\bQ$-linear combination of products $E_{2l}(\tau)I^{\rm A}(m_1,\ldots,m_{n-1};\tau)$, for $l \geq 0$ and $m_1,\ldots,m_{n-1}\geq 0$. Using the differential equation for iterated Eisenstein integrals \eqref{eqn:ieidiffeq} and the induction hypothesis, one sees that $I^{\rm A}(k_1,\ldots,k_n;\tau)$ is a $\cZ[2\pi i]$-linear combination of iterated Eisenstein integrals, plus a constant of integration, which is given by $\lim_{\tau \to i\infty}I^A(k_1,\ldots,k_n;\tau) \in \cZ[2\pi i]$, by Theorem \ref{thm:constterm}.
In conclusion, we have a unique representation
\begin{equation}
I^{\rm A}(k_1,\ldots,k_n;\tau)=\sum_{\ul{k}'} \alpha_{\ul{k}'}\cE(\ul{k}';\tau),
\end{equation}
for a finite number of multi-indices $\ul{k}'=(k_1',\ldots,k_n') \in (\bZ_{\geq 0})^n$, and $\alpha_{\ul{k}'} \in \cZ[2\pi i]$. Applying the isomorphism of Theorem \ref{thm:ieiisom} in the case $K=\cZ[2\pi i]$, we get the result.
\end{prf}
%Note that the proof shows that $\psi^A$ is the restriction of the isomorphism $\psi: \bQ\langle\cE\rangle \otimes_{\bQ} \cZ[2\pi i] \rightarrow T(\mathbf{E})^{\vee} \otimes_{\bQ} \cZ[2\pi i]$ of Theorem \ref{thm:ieiisom}.
\begin{rmk}
Recall from Section \ref{ssec:2.2} that $T(\mathbf{E})^{\vee}$ is a Hopf algebra,  whose coproduct $\Delta^{\vee}$ is given by deconcatenation. By base-extension, $\Delta^{\vee}$ naturally defines a coproduct on $T(\mathbf{E})^{\vee} \otimes_{\bQ}\cZ[2\pi i]$, which restricts to a coaction
\begin{equation}
\eZ^{\rm A} \rightarrow (T(\bfE)^{\vee} \otimes_{\bQ} \cZ[2\pi i]) \otimes_{\bQ} \eZ^{\rm A}.
\end{equation}
This coaction on $\eZ^{\rm A}$ can be seen as an elliptic analogue of the motivic coaction for motivic multiple zeta values $\cZ^{\mathfrak m}$ \cite{mtm,GonGalSym}. In fact, it is known that under a suitable isomorphism
$
\phi: \cZ^{\mathfrak m} \stackrel{\cong}\rightarrow T(\bfF)^{\vee} \otimes_{\bQ}\bQ[f_2]
$
(cf. Section \ref{ssec:1.2}),
the motivic coproduct on the Hopf algebra $\cZ^{\mathfrak m}/\zeta^{\mathfrak m}(2)$ corresponds precisely to the deconcatenation coproduct on $T(\bfF)^{\vee}$ (cf. \cite{Brown:decomp}, Section 3).
%If one could show that $\cZ^{\mathfrak m} \cong \cZ$, where $\cZ$ denotes the $\bQ$-algebra generated by the multiple zeta values, then one could define a ``full'' coproduct on $\eZ^{\rm A}$ as the tensor product
%\begin{equation}
%\Delta^{\rm geom} \otimes_{\bQ}\Delta^{\zeta^{\mathfrak m}}: \eZ^{\rm A} \rightarrow \eZ^{\rm A} \otimes_{\bQ}\eZ^{\rm A}.
%\end{equation}
\end{rmk}

\subsection{Examples} \label{ssec:4.3}
We describe some explicit examples of the decomposition map in low lengths. The case of length one is clear from Example \ref{exmp:lengthone}: we have 
\begin{equation}
\psi^{\rm A}(I^{\rm A}(2k;\tau))=\frac{2\pi iB_{2k}}{(2k)!}, \quad \psi^{\rm A}(I^{\rm A}(2k+1;\tau))=0.
\end{equation}
In what follows, we will set $\gamma_{k_1,\ldots,k_n}=\lim_{\tau \to i\infty}I^{\rm A}(k_1,\ldots,k_n;\tau)$.

\noindent
\textit{Length two:}
It follows from the differential equation \eqref{eqn:tauder} together with \eqref{eqn:lengthone} that
\begin{align} \label{eqn:lengthtwo}
I^{\rm A}(k_1,k_2;\tau)=\gamma_{k_1,k_2}&-\beta_{k_1+1,k_2}\cE(k_1+1;\tau)\notag\\
&+\beta_{k_2+1,k_1}\cE(k_2+1;\tau)\notag\\
&-(-1)^{k_2}\beta_{k_1+k_2+1,0}\cE(k_1+k_2+1;\tau)\\
&+\sum_{m=0}^{k_1+1}\binom{k_2+m-1}{m}\beta_{k_1-m+1,m+k_2}\cE(k_1-m+1;\tau)\notag\\
&-\sum_{m=0}^{k_2+1}\binom{k_1+m-1}{m}\beta_{k_2-m+1,m+k_1}\cE(k_2-m+1;\tau)\notag,
\end{align}
where $\beta_{i,j}=\alpha_i\frac{2\pi iB_j}{j!}$ if $j$ is even and $\beta_{i,j}=0$ if $j$ is odd (recall that $\alpha_i$ was defined in Theorem \ref{thm:diffeq}). In addition, comparing coefficients on both sides of \eqref{eqn:constterm}, we get
\begin{equation}
\gamma_{k_1,k_2}=\begin{cases} \displaystyle\frac{(-1)^{k_2}(2\pi i)^2}{2}\frac{B_{k_1}B_{k_2}}{k_1!k_2!}& \mbox{if $k_1 \neq 1$ or $k_2 \neq 1$,}\\\\
%\displaystyle \frac{(-1)^{k_1}}{2}\frac{B_{k_1}B_{k_2}}{k_1!k_2!} & \mbox{if $k_1=1$ and $k_2 \neq 1$}\\\\
0 & \mbox{if $k_1=k_2=1$.} \end{cases}
%\gamma_{k_1,k_2}=\begin{cases} \displaystyle\frac{B_{k_1}B_{k_2}}{2\cdot k_1!k_2!} & \mbox{if $k_2 \neq 1$}\\\\ \displaystyle-\frac{B_{k_1}B_1}{2\cdot k_1!} &  \mbox{if $k_2=1$, $k_1 \neq 1$}\\\\
%\displaystyle 0 & \mbox{if $k_1=k_2=1$.}\end{cases}
\end{equation}
One now obtains $\psi^{\rm A}(I^{\rm A}(k_1,k_2;\tau))$ by replacing in \eqref{eqn:lengthtwo} every $\cE(2m,2n;\tau)$ by $\be_{2m}^{\vee}\be_{2n}^{\vee}$ (recall that $\cE(m,n;\tau)=0$, if $m$ or $n$ is odd).

Note that, since there are no Eisenstein series of odd weight, and also since $\beta_{i,j}=0$, if $j$ is odd, we see that $I^{\rm A}(k_1,k_2;\tau)=\gamma_{k_1,k_2} \in \bQ\cdot (2\pi i)^2$, if $k_1+k_2$ is even. In particular, $I^{\rm A}(k_1,k_2;\tau)$ is, up to a power of $2\pi i$, a rational multiple of an A-elliptic multiple zeta value of length one. This is a special case of the ``length-parity theorem'' for elliptic multiple zeta values (cf. \cite{BMS}, Appendix A.1).

\noindent
\textit{Length three:}
Instead of giving a closed formula, which would be cumbersome to write down, we give a few typical examples.

Consider the A-elliptic multiple zeta value $I^{\rm A}(2,0,0;\tau)$. From \eqref{eqn:tauder}, we get
\begin{equation}
\frac{1}{2\pi i}\frac{\dd}{\dd \tau}I^{\rm A}(2,0,0;\tau)=2I^{\rm A}(3,0;\tau).
\end{equation}
On the other hand, by \eqref{eqn:constterm}, the constant term $\gamma_{2,0,0}=(2\pi i)^3\frac{B_2}{2!3!}=\frac{(2\pi i)^3}{72}$. Since $I^{\rm A}(3,0;\tau)=-2\pi i\Bigg(\cE(4;\tau)+\frac{1}{240}\cE(0;\tau)\Bigg)$, by \eqref{eqn:lengthtwo}, it follows that
\begin{equation}
I^{\rm A}(2,0,0;\tau)=2\pi i\Bigg(\frac{(2\pi i)^2}{72}-2\cE(0,4;\tau)-\frac{1}{120}\cE(0,0;\tau)\Bigg).
\end{equation}
Thus, we have
\begin{equation}
\psi^{\rm A}(I^{\rm A}(2,0,0;\tau))=2\pi i\Bigg(\frac{(2\pi i)^2}{72}-2\be^{\vee}_0\be^{\vee}_4-\frac{1}{120}\be^{\vee}_0\be^{\vee}_0\Bigg).
\end{equation}
Similarly, one shows that
\begin{align}
I^{\rm A}(0,2,0;\tau)&=2\pi i\Bigg(\frac{(2\pi i)^2}{72}+4\cE(0,4;\tau)+\frac{1}{60}\cE(0,0;\tau)\Bigg) \label{eqn:020},\\
I^{\rm A}(0,0,2;\tau)&=2\pi i\Bigg(\frac{(2\pi i)^2}{72}-2\cE(0,4;\tau)-\frac{1}{120}\cE(0,0;\tau)\Bigg) \label{eqn:002},
\end{align}
so that
\begin{align}
\psi^{\rm A}(I^{\rm A}(0,2,0;\tau))&=2\pi i\Bigg(\frac{(2\pi i)^2}{72}+4\be^{\vee}_0\be^{\vee}_4+\frac{1}{60}\be^{\vee}_0\be^{\vee}_0\Bigg),\\
\psi^{\rm A}(I^{\rm A}(0,0,2;\tau))&=2\pi i\Bigg(\frac{(2\pi i)^2}{72}-2\be^{\vee}_0\be^{\vee}_4-\frac{1}{120}\be^{\vee}_0\be^{\vee}_0\Bigg).
\end{align}
Note that $I^{\rm A}(2,0,0;\tau)=I^{\rm A}(0,0,2;\tau)$, which is an example of the reflection relations between elliptic multiple zeta values \cite{BMS,edzv}.

\textit{Length four:}
We end this section with the decomposition of $I^{\rm A}(0,1,0,0;\tau)$. This is the smallest example in which a non-trivial multiple zeta value occurs as a coefficient. Using the same procedure as before, we have by \eqref{eqn:tauder}
\begin{equation}
\frac{1}{2\pi i}\frac{\dd}{\dd \tau}I^{\rm A}(0,1,0,0;\tau)=I^{\rm A}(0,2,0;\tau)-I^{\rm A}(0,0,2;\tau),
\end{equation}
and $\gamma_{0,1,0,0}=-6\pi i\zeta(3)$ by \eqref{eqn:constterm}. Using \eqref{eqn:020} and \eqref{eqn:002}, we then get
\begin{equation}
I^{\rm A}(0,1,0,0;\tau)=2\pi i\Bigg(-3\zeta(3)+6\cE(0,0,4;\tau)+\frac{1}{40}\cE(0,0,0;\tau)\Bigg),
\end{equation}
which yields
\begin{equation}
\psi^{\rm A}(I^{\rm A}(0,1,0,0;\tau))=2\pi i\Bigg(-3\zeta(3)+6\be^{\vee}_0\be^{\vee}_0\be^{\vee}_4+\frac{1}{40}\be^{\vee}_0\be^{\vee}_0\be^{\vee}_0\Bigg).
\end{equation}
\section{The image of the decomposition map} \label{sec:5}

In the last section, we have constructed an embedding
\begin{equation}
\psi^{\rm A}: \eZ^{\rm A} \hookrightarrow T(\mathbf{E})^{\vee} \otimes_{\bQ} \cZ[2\pi i]
\end{equation}
by rewriting A-elliptic multiple zeta values as iterated Eisenstein integrals. In this section, we will see that $\psi^{\rm A}$ is not surjective, and that its image lies in a subspace associated to a certain Lie algebra of derivations $\fu$ \cite{Pol,Tsu} (see below). The key to establish this result is the differential equation for the generating series $\underline{A}(\tau)$ of A-elliptic multiple zeta values \cite{Eass}.

\subsection{A Lie algebra of derivations} \label{ssec:5.1}
Let $\cL$ be the free $\bC$-Lie algebra on the set $\{x,y\}$ (cf. e.g. \cite{SerLie}, Chapter IV). For every $k \geq 0$, define a derivation $\varepsilon_{2k}: \cL \rightarrow \cL$ by the formula
\begin{equation}
\varepsilon_{2k}(x)=\ad^{2k}(x)(y), \quad \varepsilon_{2k}(y)=\sum_{0 \leq j <k}(-1)^j[\ad^j(x)(y),\ad^{2k-1-j}(x)(y)],
\end{equation} 
where $\ad^n(x)(y):=\underbrace{[x,[x,\ldots,[x}_n,y],\ldots,]]$. Note that $\varepsilon_0$ is simply the derivation $y\frac{\partial}{\partial x} \in \mathfrak{sl}_2$, while $\varepsilon_2=-\ad([x,y])$ is an inner derivation. Also, for every $k$, we have $\varepsilon_{2k}([x,y])=0$ (cf. e.g. \cite{Pol}).

Let $\Der^0(\cL)$ be the Lie $\bQ$-algebra of derivations of $\cL$, which annihilate $[x,y]$, and let
\begin{equation}
\fu=\Lie(\varepsilon_{2k},\, k\geq 0) \subset \Der^0(\cL)
\end{equation}
be the Lie subalgebra of $\Der^0(\cL)$ generated by the $\varepsilon_{2k}$. The Lie algebra $\fu$ was first studied by Tsunogai \cite{Tsu} in a slightly different context (Galois representations of once-punctured elliptic curves).
In his master thesis \cite{Pol}, Pollack showed that relations between commutators of $\varepsilon_{2k}$'s can be traced back to modular forms for $\SL_2(\bZ)$. In particular, $\fu$ is not freely generated by the $\varepsilon_{2k}$. An equivalent formulation of this fact goes as follows. Let $U(\fu)$ be the universal enveloping algebra of $\fu$, and recall the definition of the tensor algebra $T(\bfE)$ (cf. Section \ref{ssec:2.2}). Since $\fu$ is generated by elements $\varepsilon_{2k}$ for $k \geq 0$, there exists a canonical surjection of Hopf $\bQ$-algebras
\begin{align}
T(\mathbf{E}) &\rightarrow U(\fu) \notag\\
\be_{2k} &\mapsto \varepsilon_{2k},
\end{align}
and the fact that $\fu$ is not freely generated by the $\varepsilon_{2k}$ means that this morphism is not injective. Equivalently, the dual morphism
\begin{equation} \label{eqn:dualinj}
\iota: U(\fu)^{\vee} \hookrightarrow T(\mathbf{E})^{\vee}
\end{equation}
is not surjective, where $U(\fu)^{\vee}$ denotes the graded dual of $U(\fu)$ (all $\varepsilon_{2k}$ have degree one).
\begin{exmp}
In $\fu$, we have for example the relation (cf. \cite{Pol}, eq. (3))
\begin{equation} \label{eqn:Eisrel}
[\varepsilon_{2k},\varepsilon_2]=0, \quad \forall k \geq 0,
\end{equation}
which follows from $\varepsilon_2=-\ad([x,y])$, and from the fact that every $\varepsilon_{2k}$ annihilates $[x,y]$. Since $[\varepsilon_{2k},\varepsilon_2]=\varepsilon_{2k}\circ \varepsilon_2-\varepsilon_2\circ \varepsilon_{2k}$, the relation \eqref{eqn:Eisrel} implies that a linear combination $\sum \lambda_{\underline{2k}}\be_{\underline{2k}}^{\vee}$ is contained in $\iota(U(\fu)^{\vee})$, only if $\lambda_{2,2k}=\lambda_{2k,2}$ for every $k \geq 0$.

A more interesting example is the relation (cf. \cite{Pol}, eq. (4))
\begin{equation} \label{eqn:IharaTakao}
[\varepsilon_{10},\varepsilon_4]-3[\varepsilon_8,\varepsilon_6]=\varepsilon_{10}\circ \varepsilon_4-\varepsilon_4\circ\varepsilon_{10}-3(\varepsilon_8\circ \varepsilon_6-\varepsilon_6\circ\varepsilon_8)=0,
\end{equation}
which essentially goes back to Ihara and Takao. Let $W \subset T(\bfE)^{\vee}$ be the four-dimensional subspace spanned by $\be_{10}^{\vee}\be_4^{\vee}$, $\be_4^{\vee}\be_{10}^{\vee}$, $\be_8^{\vee}\be_6^{\vee}$ and $\be_6^{\vee}\be_8^{\vee}$. Then the intersection $\iota(U(\fu)^{\vee}) \cap W$ is contained in the (three-dimensional) annihilator of \eqref{eqn:IharaTakao}, viewed as the column vector $(1,-1,-3,3)^t \in \bQ^4$. Explicitly
\begin{equation}
%\iota(U(\fu)^{\vee}) \cap W \subseteq \spn_{\bQ}\{\be_{10}^{\vee}\be_4^{\vee}+\be_4^{\vee}\be_{10}^{\vee}, \be_8^{\vee}\be_6^{\vee}+\be_6^{\vee}\be_8^{\vee},3\be_{10}^{\vee}\be_4^{\vee}+\be_8^{\vee}\be_6^{\vee} \},
\iota(U(\fu)^{\vee}) \cap W \subseteq \spn_{\bQ}\{\be_{10}^{\vee}\shuffle \be_4^{\vee}, \be_8^{\vee}\shuffle\be_6^{\vee},3\be_{10}^{\vee}\be_4^{\vee}+\be_8^{\vee}\be_6^{\vee} \},
\end{equation}
and one can show that equality holds.
\end{exmp}

\subsection{The differential equation revisited} \label{ssec:5.2}
As was already mentioned at the end of Section \ref{sec:3}, the differential equation for A-elliptic multiple zeta values can be reformulated as a differential equation for its generating series $\underline{A}(\tau)$. The precise result is the following\footnote{See also the footnote on page 7.}
%\begin{equation}
%\underline{A}(\tau)=\sum_{n \geq 0}(-1)^n\sum_{k_1,\ldots,k_n\geq 0}I^{\rm A}(k_1,\ldots,k_n;\tau)\ad^{k_n}(a)(b)\ldots\ad^{k_1}(a)(b).
%\end{equation}
%where the variables $a$ and $b$ are given by $a=2\pi ix$ and $b=y$.
\begin{thm}[\cite{Eass}, Proposition 6.2, and \cite{Eemzv}, eq. (7)] \label{thm:diffeqass}
The series $\underline{A}(\tau)$ satisfies the differential equation
\begin{equation} \label{eqn:diffeqass}
\frac{1}{2\pi i}\frac{\dd }{\dd \tau}A(\tau)=\Bigg(-\sum_{k \geq 0}E_{2k}(\tau)\widetilde{\varepsilon}_{2k}\Bigg)\Big(\underline{A}(\tau)\Big),
\end{equation}
where $E_{2k}(\tau)$ denotes the Eisenstein series \eqref{eqn:Eisen}, and \begin{equation}
\widetilde{\varepsilon}_{2k}=\begin{cases}\frac{2}{(2k-2)!}\varepsilon_{2k} & \mbox{if $k>0$,}\\ -\varepsilon_0 & \mbox{if $k=0$.}\end{cases}
\end{equation}
\end{thm}
As shown in \cite{Eemzv}, Section 4, Theorem \ref{thm:diffeqass} is equivalent to Theorem \ref{thm:diffeq}. Solving \eqref{eqn:diffeqass} iteratively, using in addition the initial condition $\underline{A}_{\infty}:=\lim_{\tau \to i\infty}\underline{A}(\tau)$, which is known explicitly by Theorem \ref{thm:constterm}, we get that
\begin{equation} \label{eqn:Aassocdiff}
\underline{A}(\tau)=g(\tau)(\underline{A}_{\infty})=g(\tau)(e^{\pi it}\Phi(\tilde{y},t)e^{2\pi i\tilde{y}}\Phi(\tilde{y},t)^{-1}),
\end{equation}
with $g(\tau)=\sum \cE(\underline{2k};\tau)\widetilde{\varepsilon}_{\underline{2k}}$, where the sum is over all multi-indices $(2k_1,\ldots,2k_n) \in (2\bZ_{\geq 0})^n$ (for all $n \geq 0$), and $\widetilde{\varepsilon}_{\underline{2k}}=\widetilde{\varepsilon}_{2k_1}\circ\ldots\circ\widetilde{\varepsilon}_{2k_n}$.

Theorem \ref{thm:diffeqass} is the key to relate A-elliptic multiple zeta values to the Lie algebra $\fu$.
\begin{thm} \label{thm:psiAEisen}
The decomposition map $\psi^{\rm A}: \eZ^{\rm A} \hookrightarrow T(\mathbf{E})^{\vee} \otimes_{\bQ}\cZ[2\pi i]$ factors through the subspace $\iota(U(\fu)^{\vee}) \otimes_{\bQ}\cZ[2\pi i]$, where $\iota: U(\fu)^{\vee} \hookrightarrow T(\mathbf{E})^{\vee}$ is the natural dual injection \eqref{eqn:dualinj}.
\end{thm}
\begin{prf}
Let $\mathcal{B}$ be a homogeneous vector space basis for $U(\fu)$, and write
\begin{equation}
g(\tau)=\sum_{b \in \mathcal{B}}\Bigg[\sum\lambda_{\underline{2k},b}\cE(\underline{2k};\tau)\Bigg] \cdot b,
\end{equation}
where $\lambda_{\underline{2k},b} \in \bQ$ and the innermost sum is finite for every $b$. Under the isomorphism $\psi: \bQ\langle\cE\rangle \rightarrow T(\mathbf{E})^{\vee}$ of Theorem \ref{thm:ieiisom}, the element $g(\tau)$ corresponds to
\begin{equation}
\psi(g(\tau))=\sum_{b \in \mathcal{B}} \Bigg[\sum\lambda_{\underline{2k},b}\be_{\underline{2k}}^{\vee}\Bigg] \cdot b,
\end{equation}
and $\psi(g(\tau))$ can be seen as a morphism
\begin{align}
U(\fu)^{\vee} &\rightarrow T(\mathbf{E})^{\vee}\notag\\
b^{\vee} &\mapsto b^{\vee}(\psi(g(\tau)))=\sum\lambda_{\underline{2k},b}\cE(\underline{2k};\tau). \label{eqn:morph}
\end{align}
This morphism is clearly dual to the natural surjection $T(\mathbf{E})\rightarrow U(\fu)$ given by $\be_{\underline{2k}} \mapsto \widetilde{\varepsilon}_{\underline{2k}}$, thus, comparing with \eqref{eqn:dualinj}, we see that the image of \eqref{eqn:morph} is equal to $\iota(U(\fu)^{\vee})$.

Now since $\underline{A}(\tau)$ is the generating series of the $I^{\rm A}(k_1,\ldots,k_n;\tau)$, the $\bQ$-span of the coefficients of $\underline{A}(\tau)=g(\tau)(\underline{A}_{\infty})$ equals $\eZ^{\rm A}$. On the other hand, by definition the image of $\psi^{\rm A}$ is equal to the $\bQ$-span of the coefficients of $\psi(g(\tau))(\underline{A}_{\infty})$, which, by the preceding discussion and the fact (cf. Theorem \ref{thm:constterm}) that the coefficients of $\underline{A}_{\infty}$ lie in $\cZ[2\pi i]$, are contained in $\iota(U(\fu)^{\vee}) \otimes_{\bQ}\cZ[2\pi i]$. Thus, the image of $\psi^{\rm A}$ is indeed contained in $\iota(U(\fu)^{\vee}) \otimes_{\bQ}\cZ[2\pi i]$.
\end{prf}
\begin{rmk}
Essentially the same result holds for the algebra $\eZ^{\rm B}$ of B-elliptic multiple zeta values. More precisely, one has a canonical embedding \cite{thesis}
\begin{equation}
\psi^{\rm B}: \eZ^{\rm B} \hookrightarrow T(\bfE)^{\vee} \otimes_{\bQ}\cZ[2\pi i].
\end{equation}

%Moreover, it is expected that $\psi^{\rm B}$ is an isomorphism \cite{LMS}.
\end{rmk}

\subsection{The Fourier subspace} \label{ssec:5.3}
In the last subsection, we have seen that
the relation between the differential equation for $\underline{A}(\tau)$ and the derivations $\widetilde{\varepsilon}_{2k}$ implies that
the image of $\psi^{\rm A}$ lies in $\iota(U(\fu)^{\vee}) \otimes_{\bQ} \cZ[2\pi i]$. In this subsection, we will see that the Fourier expansion of A-elliptic multiple zeta values (cf. Proposition \ref{prop:Fourier}) further constrains the image of $\psi^{\rm A}$.\footnote{We should note that this additional constraint is a particular feature of A-elliptic multiple zeta values. More precisely, the analogue of Theorem \ref{thm:Fou} for B-elliptic multiple zeta values is false, since B-elliptic multiple zeta values in general do not have a Fourier expansion.}
\begin{dfn} \label{dfn:Fou}
The \textit{Fourier subspace} $\bQ\langle\cE\rangle_{\rm Fou} \subset \bQ\langle\cE\rangle$ is the $\bQ$-linear subspace defined by
\begin{equation}
\bQ\langle\cE\rangle_{\rm Fou}:=\spn_{\bQ}\{\cE^0(2k_1,\ldots,2k_n;\tau)\, \vert \, n \geq 0, \, k_i \geq 0 \},
\end{equation}
where $\cE^0(2k_1,\ldots,2k_{n-1},0;\tau):=0$ and for $k_n \neq 0$, we set
\begin{equation}
\cE^0(2k_1,\ldots,2k_n;\tau):=\cE(2k_1,\ldots,2k_n;\tau)-\frac{B_{2k_n}}{4k_n}\cE(2k_1,\ldots,2k_{n-1},0;\tau).
\end{equation}
\end{dfn}
We will denote by $T(\bfE)^{\vee}_{\rm Fou}$ the subspace of $T(\bfE)^{\vee}$, which is the image of $\bQ\langle\cE\rangle_{\rm Fou}$ under the isomorphism $\psi: \bQ\langle\cE\rangle \stackrel{\cong}\rightarrow T(\bfE)^{\vee}$ of Theorem \ref{thm:ieiisom}. Note that $T(\bfE)^{\vee}_{\rm Fou}$ is a $\bQ$-subalgebra of $T(\bfE)^{\vee}$ and a left comodule under $T(\bfE)^{\vee}$,
%Note that $T(\bfE)^{\vee}_{\rm Fou}$ is a left co-ideal of $T(\bfE)^{\vee}$
i.e. the coproduct $\Delta^{\vee}$ on $T(\bfE)^{\vee}$ restricts to a morphism
\begin{equation}
T(\bfE)^{\vee}_{\rm Fou} \rightarrow T(\bfE)^{\vee} \otimes_{\bQ}T(\bfE)^{\vee}_{\rm Fou}.
\end{equation}
\begin{rmk}
The name ``Fourier subspace'' is motivated by the fact that a $\bQ$-linear combination of iterated Eisenstein integrals $\cE(2k_1,\ldots,2k_n;\tau)$ has a Fourier expansion in $q=e^{2\pi i\tau}$, if and only if it is contained in $\bQ\langle\cE\rangle_{\rm Fou}$. This follows easily from $E_{2k}(\tau)=-\frac{B_{2k}}{4k}+O(q)$, valid for $k>0$, which together with $E_0(\tau):=-1$ implies that $\cE^0(2k_1,\ldots,2k_n;\tau) \in O(q)$ (since the ideal $q\cdot \bC[[q]] \subset \bC[[q]]$ is closed under integration with respect to the measure $2\pi i\dd\tau=\dd(\log q)$).
\end{rmk}
\begin{thm} \label{thm:Fou}
The morphism $\psi^A$ maps $\eZ^{\rm A}$ into the Fourier subspace; more precisely
\begin{equation} \label{eqn:thmFou}
\psi^A: \eZ^{\rm A} \hookrightarrow \iota(U(\fu)^{\vee})_{\rm Fou} \otimes_{\bQ}\cZ[2\pi i],
\end{equation}
where $\iota(U(\fu)^{\vee})_{\rm Fou}:=\iota(U(\fu)^{\vee}) \cap T(\bfE)^{\vee}_{\rm Fou}$, and $\iota: U(\fu)^{\vee} \hookrightarrow T(\bfE)^{\vee}$ is the natural injection \eqref{eqn:dualinj}.
\end{thm}
\begin{prf}
We can rewrite $g(\tau)$ using the $\cE^0(2k_1,\ldots,2k_n;\tau)$ as follows:
\small{
\begin{alignat}{3}
&g(\tau)&&=\sum\cE(\underline{2k};\tau)\widetilde{\varepsilon}_{\underline{2k}}&&\notag\\
& &&=\sum_{k_n \neq 0}\Big(\cE^0(\underline{2k};\tau)&&+\frac{B_{2k_n}}{4k_n}\cE(2k_1,\ldots,2k_{n-1},0;\tau)\Big)\widetilde{\varepsilon}_{\underline{2k}}\notag\\
& && &&+\sum_{k_n=0}\cE(2k_1,\ldots,2k_{n-1},0;\tau)\widetilde{\varepsilon}_{2k_1}\circ \ldots \circ \widetilde{\varepsilon}_{2k_{n-1}}\circ \widetilde{\varepsilon}_0\\
& &&=\sum\cE^0(\underline{2k};\tau)\widetilde{\varepsilon}_{\underline{2k}}&&+\sum\cE(2k_1,\ldots,2k_{n-1},0;\tau)\widetilde{\varepsilon}_{2k_1} \circ \ldots \circ \widetilde{\varepsilon}_{2k_{n-1}} \circ \Big(\underbrace{\widetilde{\varepsilon}_0+\sum_{k_n \geq 1}\frac{B_{2k_n}}{4k_n}\widetilde{\varepsilon}_{2k_n}}_{=:D}\Big),\notag
\end{alignat}
}
where all sums are over the multi-indices $\underline{2k}=(2k_1,\ldots,2k_n) \in (2\bZ_{\geq 0})^n$, for all $n \geq 0$.
It is shown in the proof of \cite{Eass}, Proposition 6.3, that $D$ is a derivation of $\bC\langle\!\langle a,b\rangle\!\rangle$ that annihilates both $\tilde{y}=-\frac{\ad(a)}{e^{\ad(a)}-1}(b)$ and $t=-[a,b]$, thus it annihilates every word in $\tilde{y}$ and $t$. Since $\underline{A}_{\infty}=e^{\pi it}\Phi(\tilde{y},t)e^{2\pi i\tilde{y}}\Phi(\tilde{y},t)^{-1}$ is a power series in $\tilde{y}$ and $t$, it follows that $D(\underline{A}_{\infty})=0$. Hence, 
\begin{equation}
\underline{A}(\tau)=g(\tau)(\underline{A}_{\infty})=\Big(\sum\cE^0(\underline{2k};\tau) \widetilde{\varepsilon}_{\underline{2k}}\Big)(\underline{A}_{\infty}),
\end{equation}
and therefore every coefficient of $\underline{A}(\tau)$ is contained in $\bQ\langle\cE\rangle_{\rm Fou} \otimes_{\bQ}\cZ[2\pi i]$. Combining this with Theorem \ref{thm:psiAEisen}, the result follows.
\end{prf}
%\begin{rmk}
%The author expects that \eqref{eqn:thmFou} is an isomorphism.
%\end{rmk}
%
%
%
%\subsection{Interpretation via universal mixed elliptic motives}

\bibliography{document_typoscorrected2}

\def\cprime{$'$}
\begin{thebibliography}{10}

\bibitem{BMMS}
J.~Broedel, C.~R. Mafra, N.~Matthes, and O.~Schlotterer.
\newblock Elliptic multiple zeta values and one-loop superstring amplitudes.
\newblock {\em J. High Energy Phys.}, (7):112, front matter+41, 2015.

\bibitem{BMS}
J.~Broedel, N.~Matthes, and O.~Schlotterer.
\newblock Relations between elliptic multiple zeta values and a special
  derivation algebra.
\newblock {\em J. Phys. A}, 49(15):155203, 49, 2016.

\bibitem{mtm}
F.~Brown.
\newblock Mixed {T}ate motives over {$\mathbb Z$}.
\newblock {\em Ann. of Math. (2)}, 175(2):949--976, 2012.

\bibitem{Brown:decomp}
F.~Brown.
\newblock On the decomposition of motivic multiple zeta values.
\newblock In {\em Galois-{T}eichm\"uller theory and arithmetic geometry},
  volume~63 of {\em Adv. Stud. Pure Math.}, pages 31--58. Math. Soc. Japan,
  Tokyo, 2012.

\bibitem{Colombia}
F.~Brown.
\newblock Iterated integrals in quantum field theory.
\newblock In {\em Geometric and topological methods for quantum field theory},
  pages 188--240. Cambridge Univ. Press, Cambridge, 2013.

\bibitem{MMV}
F.~Brown.
\newblock Multiple modular values for ${\SL_2(\bZ)}$.
\newblock arXiv:1407.5167, 2014.

\bibitem{BL}
F.~Brown and A.~Levin.
\newblock Multiple elliptic polylogarithms.
\newblock arXiv:1110.6917, 2011.

\bibitem{CEE}
D.~Calaque, B.~Enriquez, and P.~Etingof.
\newblock Universal {KZB} equations: the elliptic case.
\newblock In {\em Algebra, arithmetic, and geometry: in honor of {Y}u. {I}.
  {M}anin. {V}ol. {I}}, volume 269 of {\em Progr. Math.}, pages 165--266.
  Birkh\"auser Boston, Inc., Boston, MA, 2009.

\bibitem{Ch}
K.~T. Chen.
\newblock Iterated path integrals.
\newblock {\em Bull. Amer. Math. Soc.}, 83(5):831--879, 1977.

\bibitem{Del}
P.~Deligne.
\newblock Le groupe fondamental de la droite projective moins trois points.
\newblock In {\em Galois groups over ${\bf Q}$ ({B}erkeley, {CA}, 1987)},
  volume~16 of {\em Math. Sci. Res. Inst. Publ.}, pages 79--297. Springer, New
  York, 1989.

\bibitem{Dr}
V.~G. Drinfel{\cprime}d.
\newblock On quasitriangular quasi-{H}opf algebras and on a group that is
  closely connected with {${\rm Gal}(\overline{\bf Q}/{\bf Q})$}.
\newblock {\em Algebra i Analiz}, 2(4):149--181, 1990.

\bibitem{Eass}
B.~Enriquez.
\newblock Elliptic associators.
\newblock {\em Selecta Math. (N.S.)}, 20(2):491--584, 2014.

\bibitem{Eemzv}
B.~Enriquez.
\newblock Analogues elliptiques des nombres multiz{\'e}tas.
\newblock {\em Bull. Soc. Math. France}, 144(3):395--427, 2016.

\bibitem{FurStab}
H.~Furusho.
\newblock The multiple zeta value algebra and the stable derivation algebra.
\newblock {\em Publ. Res. Inst. Math. Sci.}, 39(4):695--720, 2003.

\bibitem{GonGalSym}
A.~B. Goncharov.
\newblock Galois symmetries of fundamental groupoids and noncommutative
  geometry.
\newblock {\em Duke Math. J.}, 128(2):209--284, 2005.

\bibitem{HM}
R.~{Hain} and M.~{Matsumoto}.
\newblock {Universal Mixed Elliptic Motives}.
\newblock arXiv:1512.03975, 2015.

\bibitem{Hai}
R.~M. Hain.
\newblock The geometry of the mixed {H}odge structure on the fundamental group.
\newblock In {\em Algebraic geometry, {B}owdoin, 1985 ({B}runswick, {M}aine,
  1985)}, volume~46 of {\em Proc. Sympos. Pure Math.}, pages 247--282. Amer.
  Math. Soc., Providence, RI, 1987.

\bibitem{LR}
A.~Levin and G.~Racinet.
\newblock Towards multiple elliptic polylogarithms.
\newblock arXiv:math/0703237, 2007.

\bibitem{LMS}
P.~Lochak, N.~Matthes, and L.~Schneps.
\newblock {Elliptic multiple zeta values and the elliptic double shuffle relations}.
\newblock arXiv:1703.09410.

\bibitem{Man}
Y.~I. Manin.
\newblock Iterated integrals of modular forms and noncommutative modular
  symbols.
\newblock In {\em Algebraic geometry and number theory}, volume 253 of {\em
  Progr. Math.}, pages 565--597. Birkh\"auser Boston, Boston, MA, 2006.

\bibitem{thesis}
N.~Matthes.
\newblock {\em {Elliptic multiple zeta values}}.
\newblock PhD thesis, Universit\"at Hamburg, 2016.

\bibitem{edzv}
N.~Matthes.
\newblock Elliptic double zeta values.
\newblock {\em J. Number Theory}, 171:227--251, 2017.

\bibitem{Pol}
A.~Pollack.
\newblock {Relations between derivations arising from modular forms}.
\newblock Master's thesis, Duke University, 2009.

\bibitem{SerLie}
J.-P. Serre.
\newblock {\em Lie algebras and {L}ie groups}, volume 1500 of {\em Lecture
  Notes in Mathematics}.
\newblock Springer-Verlag, Berlin, 2006.
\newblock 1964 lectures given at Harvard University, Corrected fifth printing
  of the second (1992) edition.

\bibitem{Tsu}
H.~Tsunogai.
\newblock On some derivations of {L}ie algebras related to {G}alois
  representations.
\newblock {\em Publ. Res. Inst. Math. Sci.}, 31(1):113--134, 1995.

\bibitem{W}
A.~Weil.
\newblock {\em Elliptic functions according to {E}isenstein and {K}ronecker}.
\newblock Springer-Verlag, Berlin-New York, 1976.
\newblock Ergebnisse der Mathematik und ihrer Grenzgebiete, Band 88.

\bibitem{Zag}
D.~Zagier.
\newblock Periods of modular forms and {J}acobi theta functions.
\newblock {\em Invent. Math.}, 104(3):449--465, 1991.

\end{thebibliography}

\end{document}